\documentclass{amsart}
\usepackage{amsmath}
\usepackage{amscd}
\usepackage{amssymb}

\newcommand{\cal}{\mathcal}
\renewcommand{\div}{\operatorname{div}}
\newcommand{\im}{\operatorname{im}}
\newcommand{\vol}{\operatorname{vol}}

\newcommand{\tr}{\operatorname{tr}}

\newcommand{\CP}{\mathbb{CP}}

\newcommand{\RR}{\mathbb{R}}

\newcommand{\DDD}{{\cal D}}

\newtheorem{theorem}{Theorem}[section]
\newtheorem{theorem/definition}{Theorem/Definition}[section]
\newtheorem{proposition}{Proposition}[section]
\newtheorem{lemma}{Lemma}[section]

\theoremstyle{remark}
\newtheorem{remark}{Remark}[section]

\theoremstyle{definition}
\newtheorem{definition}{Definition}[section]
\newtheorem{example}{Example}[section]

\begin{document}
\title
{Recent Progress on Ricci Solitons}
\author{Huai-Dong Cao}
\address{Department of Mathematics\\ Lehigh University\\
Bethlehem, PA 18015} \email{huc2@lehigh.edu}

\begin{abstract}
In recent years, there has seen much interest and increased
research activities in Ricci solitons. Ricci solitons are natural
generalizations of Einstein metrics. They are also special
solutions to Hamilton's Ricci flow and play important roles in the
singularity study of the Ricci flow. In this paper, we survey some
of the recent progress on Ricci solitons.

\end{abstract}

\maketitle
\date{}


\footnotetext[1]{Research partially supported by NSF grants
DMS-0354621 and DMS-0506084.}

The concept of {\it Ricci solitons} was introduced by Hamilton
\cite{Ha88} in mid 80's. They are natural generalizations of
Einstein metrics. Ricci solitons also correspond to self-similar
solutions of Hamilton's Ricci flow \cite{Ha82} and often arise as
limits of dilations of singularities in the Ricci flow
\cite{Ha93E, Cao94, CZ00, Se}. They can be viewed as fixed points
of the Ricci flow, as a dynamical system, on the space of
Riemannian metrics modulo diffeomorphisms and scalings. Ricci
solitons are of interests to physicists as well and are called
{\it quasi-Einstein} metrics in physics literature (see, e.g.,
\cite{Fr}). In this paper, we survey some of the recent progress
on Ricci solitons as well as the role they play in the singularity
study of the Ricci flow. This paper can be regarded as an update
of the article \cite{Cao05} written by the author a few years ago.

\section{Ricci Solitons}

\subsection{Ricci Solitons} Recall that a Riemannian metric $g_{ij}$ is {\it
Einstein} if its Ricci tensor $R_{ij}=\rho g_{ij}$ for some
constant $\rho$. A smooth $n$-dimensional manifold $M^n$ with an
Einstein metric $g$ is an {\it Einstein manifold}. Ricci solitons,
introduced by Hamilton \cite{Ha88}, are natural generalizations of
Einstein metrics.

\begin{definition}
A complete Riemannian metric $g_{ij}$ on a smooth manifold $M^n$
is called a {\it Ricci soliton} if there exists a smooth vector
field $V=(V^{i})$ such that the Ricci tensor $R_{ij}$ of the
metric $g_{ij}$ satisfies the equation
$$R_{ij}+ \frac{1}{2} (\nabla_iV_j +\nabla_jV_i)=\rho g_{ij},
\eqno(1.1)$$
for some constant $\rho$.
Moreover, if $V$ is a gradient vector field, then we have a {\it
gradient Ricci soliton}, satisfying the equation
$$R_{ij}+\nabla_i\nabla_jf=\rho g_{ij},\eqno(1.2)$$
for some smooth function $f$ on $M$.
For $\rho=0$ the Ricci soliton is {\it steady}, for $\rho>0$ it is
{\it shrinking} and for $\rho<0$ {\it expanding}. The function $f$
is called a {\it potential function} of the Ricci soliton.

\end{definition}

Since $\nabla_iV_j +\nabla_jV_i$ is the Lie derivative $L_Vg_{ij}$
of the metric $g$ in the direction of $V$, we also write the Ricci
soliton equations (1.1) and (1.2) as
$$Rc+\frac{1}{2} L_Vg=\rho g \quad\mbox{and} \quad Rc+\nabla^{2}f=\rho
g \eqno(1.3)$$ respectively.

When the underlying manifold is a complex manifold, we have the
corresponding notion of K\"ahler-Ricci solitons.
\begin{definition}
A complete K\"ahler metric $g_{\alpha\bar\beta}$ on a complex
manifold $X^n$ of complex dimension $n$ is called a {\it
K\"ahler-Ricci soliton} if there exists a holomorphic vector field
$V=(V^{\alpha})$ on $X$ such that the Ricci tensor
$R_{\alpha\bar\beta}$ of the metric $g_{\alpha\bar\beta}$
satisfies the equation
$$R_{\alpha\bar\beta}+\frac{1}{2}(\nabla_{\bar\beta}V_{\alpha}+
\nabla_{\alpha} V_{\bar\beta})=\rho g_{\alpha\bar\beta}
\eqno(1.4)$$
for some (real) constant $\rho$. It is called a gradient
K\"ahler-Ricci soliton if the holomorphic vector field $V$ comes
from the gradient vector field of a real-valued function $f$ on
$X^n$ so that
$$R_{\alpha\bar\beta}+\nabla_{\alpha}\nabla_{\bar\beta} f=\rho
g_{\alpha\bar\beta}, \quad \mbox{and} \quad
\nabla_{\alpha}\nabla_{\beta} f=0. \eqno(1.5)$$
 Again, for
$\rho=0$ the soliton is {\it steady}, for $\rho>0$ it is shrinking
and for $\rho<0$ expanding.
\end{definition}

Note that the case $V=0$ (i.e., $f$ being a constant function) is
an Einstein (or K\"ahler-Einstein) metric. Thus Ricci solitons are
natural extensions of Einstein metrics. In fact, we will see below
that there are no non-Einstein compact steady or expanding Ricci
solitons. Also, by a suitable scale of the metric $g$, we can
normalize $\rho=0, +1/2$, or $-1/2$.

\begin{lemma} {\bf (Hamilton \cite{Ha95F})} Let $g_{ij}$ be a complete
gradient Ricci soliton with potential function $f$. Then we have
$$R+|\nabla f|^2-2\rho f=C \eqno(1.6)$$ for some constant $C$. Here $R$
denotes the scalar curvature.
\end{lemma}

\begin{proof}
Let $g_{ij}$ be a complete gradient Ricci soliton on a manifold
$M^n$ so that there exists a potential function $f$ such that the
soliton equation (1.2) holds. Taking the covariant derivatives and
using the commutating formula for covariant derivatives, we obtain
$$\nabla_iR_{jk}-\nabla_jR_{ik}+R_{ijkl}\nabla_lf=0.$$
Taking the trace on $j$ and $k$, and using the contracted second
Bianchi identity
$$\nabla_j R_{ij}=\frac {1}{2}\nabla_i R,$$
we get
$$\nabla_iR=2R_{ij}\nabla_jf. \eqno(1.7)$$
Thus
$$\nabla_i(R+|\nabla f|^2-2\rho f)=2(R_{ij}+\nabla_i\nabla_jf-\rho g_{ij})\nabla_jf=0.$$
Therefore $$R+|\nabla f|^2-2\rho f=C$$ for some constant $C$.
\end{proof}

\begin{proposition}  {\bf (cf. Hamilton \cite{Ha95F}, Ivey \cite{Iv})} On a compact manifold
$M^n$, a gradient steady or expanding Ricci soliton is necessarily
an Einstein metric.
\end{proposition}

\begin{proof}
Taking the trace in equation (1.2), we get
$$R+\Delta f =n\rho. \eqno(1.8)$$

Taking the difference of (1.6) in Lemma 1.1 and (1.8), we get
$$\Delta f -|\nabla f|^2+2\rho f=n\rho-C.$$
When $M$ is compact and $\rho\leq 0$, it follows from the maximum
principle that $f$ must be a constant and hence $g_{ij}$ is a
Einstein metric.
\end{proof}
More generally, we have

\begin{proposition}   Any compact steady or expanding Ricci soliton must be
Einstein.
\end{proposition}

\begin{proof} This follows from Proposition 1.1 and Perelman's results that
any compact Ricci soliton is necessarily a gradient soliton (see
Propositions 2.1-2.4).
\end{proof}

For compact shrinking Ricci solitons in low dimensions, we have

\begin{proposition} {\bf (Hamiton \cite{Ha88}\footnote{See alternative proofs in (Proposition 5.21,
\cite{CK}) or (Proposition 5.1.10, \cite{CZ05}), and \cite{CLT}.}
for $n=2$, Ivey \cite{Iv}\footnote{See \cite{ELM} for alternative
proofs} for $n=3$)} In dimension $n\le 3$, there are no compact
shrinking Ricci solitons other than those of constant positive
curvature.
\end{proposition}

\subsection{Examples of Ricci Solitons}

When $n\ge 4$, there exist nontrivial compact gradient shrinking
solitons. Also, there exist complete noncompact Ricci solitons
(steady, shrinking and expanding) that are not Einstein. Below we
list a number of such examples. It turns out most of the examples
are rotationally symmetric and gradient, and all the known
examples of nontrivial {\it shrinking} solitons so far are
K\"ahler.

\begin{example} {\bf (Compact gradient K\"ahler shrinkers)}
For real dimension $4$, the first example of a compact shrinking
soliton was constructed in early 90's by Koiso \cite{Ko} and the
author \cite{Cao94}\footnote{The author's construction was carried
out in 1991 at Columbia University. When he told his construction
to S. Bando that year in New York, he also learned the work of
Koiso from Bando.} on compact complex surface $\Bbb
 CP^2\#(-\Bbb CP^2)$, where $(-\Bbb CP^2)$ denotes the complex
projective space with the opposite orientation. This is a gradient
K\"ahler-Ricci soliton, has $U(2)$ symmetry and positive Ricci
curvature. More generally, they found $U(n)$-invariant
K\"ahler-Ricci solitons on twisted projective line bundle over
$\Bbb CP^{n-1}$ for $n\geq 2$.

\end{example}

\begin{remark} If a compact K\"ahler manifold $M$ admits a non-trivial
K\"ahler shrinker then $M$ is Fano (i.e., the first Chern class
$c_1(M)$ of $M$ is positive), and the Futaki-invariant \cite{Fu}
is nonzero.
\end{remark}

\begin{example} {\bf (Compact toric gradient K\"ahler shrinkers)}
In \cite{WZ}, Wang-Zhu found a gradient K\"ahler-Ricci soliton on
$\Bbb CP^2\#2(-\Bbb CP^2)$ which has $U(1)\times U(1)$ symmetry.
More generally, they proved the existence of gradient
K\"ahler-Ricci solitons on all Fano toric varieties of complex
dimension $n\geq 2$ with non-vanishing Futaki invariant.
\end{example}

\begin{example} {\bf (Noncompact gradient K\"ahler shrinkers)}
Feldman-Ilmanen-Knopf \cite{FIK} found the first complete
noncompact $U(n)$-invariant shrinking gradient K\"ahler-Ricci
solitons, which are cone-like at infinity. It has positive scalar
curvature but the Ricci curvature does not have a fixed sign.
\end{example}

\begin{example}{\bf (The cigar soliton)}
In dimension two, Hamilton \cite{Ha88} discovered the first
example of a complete noncompact steady soliton on $\mathbb R^2$,
called the {\it cigar soliton}, where the metric is given by
$$ ds^2=\frac {dx^2 +dy^2} {1+x^2+y^2}$$ with potential function
$$f=-\log (1+x^2+y^2).$$
The cigar has positive (Gaussian) curvature and linear volume
growth, and is asymptotic to a cylinder of finite circumference at
$\infty$.
\end{example}

\begin{example}{\bf (The Bryant soliton)}
In the Riemannian case, higher dimensional examples of noncompact
gradient steady solitons were found by Robert Bryant \cite{BR} on
$\Bbb R^n$ ($n\geq 3$). They are rotationally symmetric and have
positive sectional curvature. Furthermore, the geodesic sphere
$S^{n-1}$ of radius $s$ has the diameter on the order $\sqrt{s}$.
Thus the volume of geodesic balls $B_r(0)$ grow on the order of
$r^{(n+1)/2}$.
\end{example}

\begin{example}{\bf (Noncompact gradient steady K\"ahler solitons)}
In the K\"ahler case, the author \cite{Cao94} found two examples
of complete rotationally noncompact gradient steady K\"ahler-Ricci
solitons

(a) On $\Bbb C^n$ (for $n=1$ it is just the cigar soliton). These
examples are $U(n)$ invariant and have positive sectional
curvature. It is interesting to point out that the geodesic sphere
$S^{2n-1}$ of radius $s$ is an $S^1$-bundle over $\Bbb CP^{n-1}$
where the diameter of $S^1$ is on the order $1$, while the
diameter of $\Bbb CP^{n-1}$ is on the order $\sqrt{s}$. Thus the
volume of geodesic balls $B_r(0)$ grow on the order of $r^{n}$,
$n$ being the complex dimension. Also, the curvature $R(x)$ decays
like $1/r$.

(b) On the blow-up of $\Bbb C^n/\Bbb Z_n$ at the origin. This is
the same space on which Eguchi-Hansen \cite{EH} ($n=2$) and Calabi
\cite{Calabi} ($n\geq 2$) constructed examples of Hyper-K\"ahler
metrics. For $n=2$, the underlying space is the canonical line
bundle over $\mathbb{CP}^1$.
\end{example}

\begin{example}{\bf (Noncompact gradient expanding K\"ahler solitons)}
In \cite{Cao97}, the author constructed a one-parameter family of
complete noncompact expanding solitons on $\Bbb C^n$. These
expanding K\"ahker-Ricci solitons all have $U(n)$ symmetry and
positive sectional curvature, and are cone-like at infinity.

More examples of complete noncompact K\"ahler-Ricci expanding
solitons were found by Feldman-Ilmanen-Knopf \cite{FIK} on
"blow-ups" of $\Bbb C^n/\Bbb Z_k$, $k=n+1, n+2, \cdots$. (See also
Pedersen et al \cite{PTV}.)
\end{example}

\begin{example}{\bf (Sol and Nil solitons)}
Non-gradient expanding Ricci solitons on Sol and Nil manifolds
were constructed by J. Lauret \cite{La} and Baird-Laurent
\cite{BL}.
\end{example}

\begin{example} {\bf (Warped products)}
Using doubly warped product and multiple warped product
constructions, Ivey \cite{Iv94} and Dancer-Wang \cite{DW2}
produced noncompact gradient steady solitons, which generalize the
construction of Bryant's soliton.  Also, Gastel-Kronz \cite{GK}
produced a two-parameter family (doubly warped product metrics) of
gradient expanding solitons on $\Bbb R^{m+1}\times N$, where
$N^{n}$ ($n\ge2$) is an Einstein manifold with positive scalar
curvature.
\end{example}

\begin{example}
Very recently, Dancer-Wang \cite{DW1} produced new examples of
gradient shrinking, steady and expanding K\"ahler solitons on
bundles over the product of Fano K\"ahler-Einstein manifolds,
generalizing those in Examples 1.1, 1.3, 1.6, 1.7 and those by
Pedersen et al \cite{PTV}.
\end{example}

We conclude our examples with
\begin{example} {\bf (Gaussian solitons)}
$(\Bbb R^n, g_0)$ with the flat Euclidean metric can be also
equipped with both shrinking and expanding gradient Ricci
solitons, called the Gaussian shrinker or expander.

(a) $(\Bbb R^n, g_0, |x|^2/4)$ is a gradient shrinker with
potential function $f=|x|^2/4$:

$$Rc+\nabla^2 f=\frac{1} {2}g_0.$$

(b) $(\Bbb R^n, g_0, -|x|^2/4)$ is a gradient expander with
potential function $f=-|x|^2/4$:
$$Rc+\nabla^2 f=-\frac{1} {2}g_0.$$
\end{example}

\begin{remark} We'll see later that the Gaussian shrinker is very special
because it has the largest reduced volume $\tilde{V}=1$ (see
Section 4.2)
\end{remark}

\section{Variational Structures}

In this section we describe Perelman's $\mathcal F$-functional and
$\mathcal W$-functional and the associated $\lambda$-energy and
$\nu$-energy respectively. The critical points of the
$\lambda$-energy (respectively $\nu$-energy) are precisely given
by compact gradient steady (respectively shrinking) solitons. We
also consider the $\mathcal {W_{-}}$-functional and the
corresponding $\nu_{-}$-energy introduced by Feldman-Ilmanen-Ni
\cite{FIN} whose critical points are expanding solitons.
Throughout this section we assume that $M^n$ is a compact smooth
manifold.

\subsection{The $F$-functional and $\lambda$-energy} In \cite{P1}
Perelman considered the functional $${\mathcal{F}}
(g_{ij},f)=\int_M(R+|\nabla f|^2)e^{-f}dV$$ defined on the space
of Riemannian metrics and smooth functions on $M$. Here $R$ is the
scalar curvature and $f$ is a smooth function on $M^n$. Note that
when $f=0$, $F$ is simply the total scalar curvature of $g$, or
the Einstein-Hilbert action on the space of Riemannian metrics on
$M$.

\begin {lemma} {\bf (First Variation Formula of $F$-functional, Perelman \cite{P1})}
If $\delta g_{ij}=v_{ij}$ and $\delta f=\phi$ are variations of
$g_{ij}$ and $f$ respectively, then the first variation of
${\mathcal{F}}$ is given by
$$\delta{\mathcal{F}}(v_{ij},\phi)=\int_M[-v_{ij}(R_{ij}+\nabla_i\nabla_jf)+(\frac{v}{2}-\phi)(2\Delta
f-|\nabla f|^2+R)]e^{-f}dV$$ where $v=g^{ij}v_{ij}$.
\end{lemma}

Next we consider the associated energy
$$
\lambda(g_{ij})=\inf\{{\mathcal{F}}(g_{ij},f): f\in C^{\infty}(M),
\int_Me^{-f}dV=1\}.
$$
Clearly  $\lambda(g_{ij})$ is invariant under diffeomorphisms. If
we set $u=e^{-f/2}$, then the functional ${\mathcal{F}}$ can be
expressed as
$${\mathcal{F}}=\int_M(R u^2+4|\nabla
u|^2)dV.$$  Thus
$$\lambda (g_{ij})=\inf\{\int_M(R u^2+4|\nabla
u|^2)dV: \int_M u^2 dV=1\},$$ the first eigenvalue of the operator
$-4\Delta +R$. Let $u_0>0$ be a first eigenfunction of the
operator $-4\Delta +R$ so that
$$ -4\Delta u_0+Ru_0=\lambda(g_{ij})u_0.$$ Then $f_0=-2\log u_0$ is a
minimizer of $\lambda (g_{ij})$:
$$\lambda(g_{ij})={\mathcal{F}}(g_{ij},f_0).$$
Note that $f_0$ satisfies the equation
$$ -2\Delta f_0+|\nabla f_0|^2-R=\lambda(g_{ij}). \eqno(2.1)$$

For any symmetric 2-tensor $h=h_{ij}$, consider the variation
$g_{ij}(s)=g_{ij}+sh_{ij}$. It is an easy consequence of Lemma 2.1
and Eq. (2.1) that the first variation $\mathcal D_g\lambda(h)$ of
$\lambda(g_{ij})$ is given by
$$\left.\frac{d}{ds}\right|_{s=0}\lambda(g_{ij}(s))= \int
-h_{ij}(R_{ij}+\nabla_i\nabla_j f)e^{-f} dV, \eqno(2.2)$$
where $f$ is a minimizer of $\lambda(g_{ij}$. In particular, {\it
the critical points of $\lambda$ are precisely steady gradient
Ricci solitons}.

Note that, by diffeomorphism invariance of $\lambda$, $\mathcal
D_g\lambda$ vanishes on any Lie derivative $h_{ij}=\frac{1}{2}
L_Vg_{ij}$, and hence on $\nabla_i\nabla_j f=\frac{1}{2}L_{\nabla
f}g_{ij}$. Thus, by inserting $h=-2(Ric+\nabla^2f)$ in Eq. (2.2)
one recovers the following result of Perelman \cite{P1}.

\begin{proposition} Suppose $g_{ij}(t)$ is a solution to the Ricci flow on
a compact manifold $M^n$. Then $\lambda(g_{ij}(t))$ is
nondecreasing in $t$ and the monotonicity is strict unless we are
on a steady gradient soliton. In particular, a (compact) steady
Ricci soliton is necessarily a gradient soliton.
\end{proposition}

We remark that by considering the quantity
$$\bar{\lambda}(g_{ij})=\lambda(g_{ij})(Vol(g_{ij}))^{\frac{2}{n}},$$
which is a scale invariant version of $\lambda(g_{ij})$, Perelman
\cite{P1} also showed the following result.

\begin{proposition} $\bar{\lambda}(g_{ij})$ is nondecreasing along the Ricci flow
whenever it is nonpositive; moreover, the monotonicity is strict
unless we are on a gradient expanding soliton. In particular, any
(compact) expanding Ricci soliton is necessarily a gradient
soliton.
\end{proposition}

\subsection{The $\mathcal{W}$-functional and $\nu$-energy}

In order to study shrinking Ricci solitons, Perelman \cite{P1}
introduced the $\mathcal{W}$-functional
$$
\mathcal{W}(g_{ij},f,\tau)=\int_M [\tau (R+|\nabla f|^2)+f-n](4\pi
\tau)^{-\frac{n}{2}}e^{-f}dV,
$$
where $g_{ij}$ is a Riemannian metric, $f$ a smooth function on
$M^n$, and $\tau$ a positive scale parameter. Clearly the
functional $\mathcal{W}$ is invariant under simultaneous scaling
of $\tau$ and $g_{ij}$ (or equivalently the parabolic scaling),
and invariant under diffeomorphism. Namely, for any positive
number $a$ and any diffeomorphism $\varphi$ we have
$$
\mathcal{W}(a\varphi^*g_{ij},\varphi^*f,a
\tau)=\mathcal{W}(g_{ij},f,\tau).
$$

\begin{lemma} {\bf (First Variation of $\mathcal{W}$-functional, Perelman \cite{P1})}
If $v_{ij}=\delta g_{ij},\;
\phi=\delta f,\;\mbox{and}\; \eta=\delta\tau$, then
$$\arraycolsep=1.5pt\begin{array}{rcl}\delta \mathcal{W}(v_{ij},\phi,\eta)&=&\int_M-\tau
v_{ij}(R_{ij}+\nabla_if\nabla_jf-\frac{1}{2\tau}g_{ij})(4\pi\tau)^{-\frac{n}{2}}e^{-f}dV\\[4mm]
&&+\int_M(\frac{v}{2}-\phi-\frac{n}{2\tau}\eta)[\tau(R+2\Delta f
-|\nabla f|^2)+f-n-1](4\pi\tau)^{-\frac{n}{2}}e^{-f}dV\\[4mm]
&&+\int_M \eta(R+|\nabla
f|^2-\frac{n}{2\tau})(4\pi\tau)^{-\frac{n}{2}}e^{-f}dV.
\end{array}
$$ Here $v=g^{ij}v_{ij}$ as before.
\end{lemma}

Similar to the $\lambda$-energy, we can consider
$$\mu(g_{ij},\tau)=\inf\{\mathcal{W}(g_{ij},f,\tau): f\in C^\infty(M),
(4\pi\tau)^{-\frac{n}{2}}\int_M e^{-f}dV=1\}. \eqno (2.3)
$$
Note that if we let $u=e^{-f/2}$, then the functional
$\mathcal{W}$ can be expressed as
$$\mathcal{W}(g_{ij},f,\tau)=\int_M[\tau(Ru^2+4|\nabla u|^2])-u^2\log
u^2-nu^2](4\pi\tau)^{-\frac{n}{2}}dV,
$$
and the constraint $\int_M(4\pi\tau)^{-\frac{n}{2}}e^{-f}dV=1$
becomes $\int_Mu^2(4\pi\tau)^{-\frac{n}{2}}dV=1$. Therefore
$\mu(g_{ij},\tau)$ corresponds to the best constant of a
logarithmic Sobolev inequality.

Since the nonquadratic term is subcritical (in view of Sobolev
exponent), it is rather straightforward to show that
$\mu(g_{ij},\tau)$ is achieved by some nonnegative function $u\in
H^1(M)$ which satisfies the Euler-Lagrange equation

$$\tau (-4\Delta u+Ru)-2u\log u-nu=\mu(g_{ij},\tau)u. \eqno(2.4)
$$
One can further show that the minimizer $u$ is positive and smooth
(see Rothaus \cite{Ro}). This is equivalent to say that
$\mu(g_{ij},\tau)$ is achieved by some minimizer $f$ satisfying
the nonlinear equation
$$\tau (2\Delta f-|\nabla f|^2+R)+f-n=\mu(g_{ij},\tau). \eqno(2.5)
$$

\begin{proposition} {\bf (Perelman \cite{P1})} Suppose $g_{ij}(t)$, $0\leq t< T$
is a solution to the Ricci flow on a compact manifold $M^n$. Then
$\mu(g_{ij}(t),T-t)$ is nondecreasing in $t$; moveover, the
monotonicity is strict unless we are on a shrinking gradient
soliton. In particular, any (compact) shrinking Ricci soliton is
necessarily a gradient soliton.
\end{proposition}

\begin{remark} Recently, Naber \cite{Na} has shown that if $(M^n,
g)$ is a complete noncompact shrinking Ricci soliton with bounded
curvature $|Rm|<C$ with respect to some smooth vector field $V$,
then there exists a smooth function $f$ on $M$ such that $(M^n,
g)$ is a gradient soliton with $f$ as a potential function. This
in particular means that $V=\nabla f +X$ for some Killing field
$X$ on $M$.
\end{remark}

The associated $\nu$-energy is defined by
$$
\nu(g_{ij})=\inf\{\mathcal W(g,f,\tau): f\in C^\infty(M), \tau>0,
(4\pi\tau)^{-\frac{n}{2}}\int e^{-f}dV=1\}.
$$
One checks that $\nu(g_{ij})$ is realized by a pair $(f,\tau)$
that solve the equations
$$\tau(-2\Delta f+|Df|^2-R)-f+n+\nu=0,\qquad
(4\pi\tau)^{-\frac{n}{2}}\int fe^{-f}=\frac{n}{2}+\nu.
\eqno(2.6)$$

Consider variations $g_{ij}(s)=g_{ij}+sh_{ij}$ as before. Using
Lemma 2.2 and (2.6), one calculates the first variation $\mathcal
D_g\nu(h)$ to be
\begin{align*}
\left.\frac{d}{ds}\right|_{s=0}\nu(g_{ij}(s))=(4\pi\tau)^{-\frac{n}{2}}\int
-h_{ij}[\tau(R_{ij}+\nabla_i\nabla_j f)-\frac{1}{2} g_{ij}]e^{-f}
dV.
\end{align*}
A stationary point of $\nu$ thus satisfies
\begin{align*}
R_{ij}+\nabla_i\nabla_j f-\frac{1}{2\tau}g_{ij}=0,
\end{align*}
which says that $g_{ij}$ is a gradient shrinking Ricci soliton.

As before, $\mathcal D_g\nu(h)$ vanishes on Lie derivatives. By
scale invariance it also vanishes on multiplies of the metric.
Inserting $h_{ij}=-2(R_{ij}+\nabla_i\nabla_j
f-\frac{1}{2\tau}g_{ij})$, one recovers Perelman's  formula that
finds that $\nu(g_{ij}(t))$ is monotone on a Ricci flow, and
constant if and only if $g_{ij}(t)$ is a gradient shrinking Ricci
soliton.

\subsection{The $\mathcal{W_{-}}$-functional and $\nu_{-}$-energy}

In \cite{FIN}, Feldman-Ilmanen-Ni introduced the dual
$\mathcal{W_{-}}$-functional (corresponding to expanders)
$$
\mathcal{W_{-}}(g_{ij},f,\sigma)=\int_M [\sigma (R+|\nabla
f|^2)-(f-n)](4\pi \sigma)^{-\frac{n}{2}}e^{-f}dV,
$$ the $\mu_{-}$-energy

$$\mu_{-}(g_{ij},\sigma)=\inf\{\mathcal{W}_{-}(g_{ij},f,\tau): f\in C^\infty(M),
(4\pi\sigma)^{-\frac{n}{2}}\int_M e^{-f}dV=1\},
$$
 and the corresponding $\nu_{-}$-entropy
$$
\nu_{-}(g_{ij})=\sup_{\sigma >0}\{\mu_{-}(g_{ij},\sigma)\}.
$$ Here, $\sigma$ is a positive parameter. They proved that

\begin{proposition} {\bf (Feldman-Ilmanen-Ni \cite{FIN})}

(a) $\mu_{-}(g_{ij},\sigma)$ is achieved by a unique $f$;
$\mu_{-}(g_{ij}(t), t-t_{0})$ is nondecreasing under the Ricci
flow; moveover, the monotonicity is strict unless we are on an
expanding gradient soliton.

(b) If $\lambda(g)<0$, then $\nu_{-}(g_{ij})$ is achieved by a
unique $\sigma$; $\nu_{-}(g_{ij}(t))$ is nondecreasing under the
Ricci flow, and is constant only on an expanding soliton.
\end{proposition}

Furthermore, if $\lambda(g)<0$ then $\nu_{-}$ is achieved by a
unique pair $(f, \sigma)$ that solve the equations
$$\sigma (-2\Delta f+|Df|^2-R)+f-n+\nu_{-}=0,\qquad
(4\pi\tau)^{-\frac{n}{2}}\int fe^{-f}=\frac{n}{2}-\nu_{-}.$$

\section{Ricci Solitons and Ricci Flow}

\subsection{Ricci solitons as self-similar solutions of the Ricci flow}

Let us first examine how Einstein metrics behave under Hamilton's
Ricci flow
$$\frac{\partial g_{ij}(t)}{\partial t}=-2R_{ij}(t).$$

If the initial metric is Ricci flat, so that $R_{ij}=0$ at $t=0$,
then clearly the metric does not change under the Ricci flow:
$g_{ij}(t)=g_{ij}(0)$. Hence any Ricci flat metric is a stationary
solution. This happens, for example, on a flat torus or on any
$K3$-surface with a Calabi-Yau metric.

If the initial metric $g_{ij}(0)$ is Einstein with positive scalar
curvature, then the metric will shrink under the Ricci flow by a
time-dependent factor. Indeed, if at $t=0$ we have
$$R_{ij}(0)=\frac{1}{2} g_{ij}(0).$$  Then
$$g_{ij}(t)=(1-t)g_{ij}(0), \eqno (3.1)$$ which shrinks homothetically to a
point as $t\rightarrow T=1$, while the scalar curvature $R\to
\infty$ like $1/(T-t)$ as $t\to T$. Note that $g(t)$ exists for
$t\in (-\infty, T)$, hence an {\sl ancient solution}.

By contrast, if the initial metric is an Einstein metric of
negative scalar curvature, the metric will expand homothetically
for all times. Suppose
$$R_{ij}(0)=- \frac{1}{2}g_{ij}(0)
$$ at $t=0$. Then the solution to the Ricci flow is given by
$$g_{ij}(t)=(1+t)g_{ij}(0).$$
Hence the evolving metric $g_{ij}(t)$ exists and expands
homothetically for all time, and the curvature will fall back to
zero like $-1/t$. Note that now the evolving metric $g_{ij}(t)$
only goes back in time to $-1$, when the metric explodes out of a
single point in a "big bang".

Now suppose we have a one-parameter group of diffeomorphisms
$\varphi_t$, $-\infty <t<\infty$, which is generated by some
vector field $V$ on $M$, and suppose $g_{ij}(t)=\varphi^*_t\hat
g_{ij}$ is a slution to the Ricci flow, called a {\it self-similar
slution}, with initial metric $\hat g_{ij}$. Then $$-2 Rc=\mathcal
L_V g$$ for all $t$. In particular, the initial metric
$g_{ij}(0)=\hat g_{ij}$ satisfies the steady Ricci soliton
equation in (1.3).

Conversely, suppose we have a steady Ricci soliton $\hat g=(\hat
g_{ij})$ on a smooth manifold $M^n$ so that
$$2 \hat Rc+\mathcal L_V \hat g=0,$$
for some smooth vector field $V=(V^i)$. {\it Assume the vector
field $V$ is complete} (i.e., $V$ generates a one-parameter group
of diffeomorphisms $\varphi_t$ of $M$). Then clearly
$$g_{ij}(t)=\varphi^*_t\hat g_{ij} \qquad -\infty<t<\infty,
$$ is a self-similar solution of the Ricci flow with $\hat
g_{ij}$ as the initial metric.

More generally, we can consider self-similar solutions to the
Ricci flow which move by diffeomorphisms and also shrinks or
expands by a (time-dependent) factor at the same time. Such
self-similar solutions correspond to either shrinking or expanding
Ricci solitons $(M, \hat g, V)$ with the vector filed $V$ being
complete. For example, a shrinking gradient Ricci soliton
satisfying the equation
\begin{align*}
\hat{R}_{ij}+\nabla_i\nabla_{j}\hat{f}-\frac{1}{2}\hat{g}_{ij}=0,
\end{align*}
with $V=\nabla \hat{f}$ complete, corresponds to the self-similar
Ricci flow solution $g_{ij}(t)$ of the form
$$g_{ij}(t):=(1-t)\varphi_t^*(g_{ij}),\qquad t<1, \eqno(3.2)$$
where $\varphi_t$ are the diffeomorphisms generated by $V/(1-t)$.
(Compare Eq. (3.2) with Eq. (3.1) for $\rho =1/2$.)

Thus, we see a complete gradient Ricci soliton with respect to
some complete vector field corresponds to the self-similar
solution of the Ricci flow it generates. For this reason we often
do not distinguish the two.

\begin{remark}
If $M^n$ is compact, then $V$ is always complete. But if $M$ is
noncompact then $V$ may not be complete in general. Recently Z.-H.
Zhang \cite{Zh1} has observed that for any {\bf complete} gradient
(steady, shrinking, or expanding) Ricci soliton $g_{ij}$with
potential function $f$, $V=\nabla f$ is a complete vector field on
$M$.

In particular, {\it a complete gradient Ricci soliton always
corresponds to the self-similar solution of the Ricci flow it
generates}
\end{remark}

\subsection{Ricci Solitons and Singularity Models of the Ricci Flow}

Ricci solitons play an important role in the study of the Ricci
flow. They are intimately related to the Li-Yau-Hamilton (also
called differential Harnack) type estimates (cf. \cite{Ha95F} and
\cite{CZ05}) in such a way that the Li-Yau-Hamilton quantity
vanishes on (expanding) Ricci solitons (see, e.g., Section 2.5 in
\cite{CZ05}). More importantly, Ricci solitons often arise as the
below-up limits of singularities in the Ricci flow which we now
describe.

Consider a solution $g_{ij}(t)$ to the Ricci flow on $M^n\times
[0,T)$, $T\le +\infty$, where either $M^n$ is compact or at each
time $t$ the metric is complete and has bounded curvature. We say
that $g_{ij}(t)$ is a {\it maximal} solution of the Ricci flow if
either $T=+\infty$ or $T<+\infty$ and the norm of its curvature
tensor $|Rm|$ is unbounded as $t\to T$. In the latter case, we say
$g_{ij}(t)$ is a singular solution to the Ricci flow.

Clearly, a round sphere $S^3$ will shrink to a point under the
Ricci flow in some finite time. Also, as Yau suggested to Hamilton
in mid 80's, if we take a dumbbell metric on $S^3$ with a neck
like $S^2\times B^1$, we expect the neck will shrink because the
positive curvature in the $S^2$ direction will dominate the
slightly negative curvature in the $B^1$ direction. In some finite
time we expect the neck will pinch off. If we dilate in space and
time at the maximal curvature point, then we expect the limit of
dilations converge to the round infinite cylinder $S^2\times
\mathbb{R}$. This intuitive picture was justified by
Angenent-Knopf \cite{AK} on $\mathbb{S}^{n+1}$ with suitable
rotationally symmetric metrics. These are examples of so called
Type I singularities. Hamilton \cite{Ha95F} also described an
intuitive picture of a degenerate neck-pinching: imagine the
dumbbell is not symmetric and one side is bigger than the other.
Then one could also pinch off a small sphere from a big one. If we
choose the sizes of the little one and the large one to be just
right, then we expect a degenerate singularity: the little sphere
pinches off and there is nothing left on the other side. H.-L Gu
and X.-P. Zhu \cite{GZ} recently verified this picture by showing
that such a degenerate neck-pinching, a Type II singularity, can
be formed on $\mathbb{S}^{n}$ with suitable rotationally symmetric
metric for all $n \geq 3$.

As in the minimal surface theory and harmonic map theory, one
usually tries to understand the structure of a singularity  by
rescaling the solution (or blow up) to obtain a sequence of
solutions and study its limit. For the Ricci flow, the theory was
first developed by Hamilton in \cite{Ha95F}.

Denote by
$$ K_{max}(t)=\sup_{x\in M}|Rm(x,t)|_{g_{ij}(t)}.
$$
According to Hamilton \cite{Ha95F}, one can classify maximal
solutions into three types; every maximal solution is clearly of
one and only one
of the following three types:\\

$\mbox{\textbf{Type I:}  } \ \ \ \ \ \ \ \ T<+\infty\ \
\mbox{and}\ \ \sup(T-t)K_{max}(t)<+\infty;\ $

\vskip 0.1cm $\mbox{\textbf{Type II(a):}} \ \ \ T<+\infty\ \
\mbox{but}\ \ \sup(T-t)K_{max}(t)=+\infty; $

\vskip 0.1cm $\mbox{\textbf{Type II(b):}} \ \ \  T=+\infty\ \
\mbox{but}\ \ \sup tK_{max}(t)=+\infty;\ \ \ \ \ \ $

\vskip 0.1cm $ \mbox{\textbf{Type III:}  }\ \ \ \ \  T=+\infty,\ \
\ \ \sup tK_{max}(t)<+\infty $\\

For each type of maximal solution, Hamilton defines a
corresponding type of limiting singularity model.

\begin{definition} A solution $g_{ij}(x,t)$ to the Ricci flow on the
manifold $M$, where either $M$ is compact or at each time $t$ the
metric $g_{ij}(\cdot,t)$ is complete and has bounded curvature, is
called a {\bf singularity model} if it is not flat and of one of
the following three types:

\vskip 0.1cm\noindent \textbf{Type I}: The solution exists for
$t\in(-\infty,\Omega)$ for some constant $\Omega$ with
$0<\Omega<+\infty$ and
$$|Rm|\leq\Omega/(\Omega-t)
$$
\hskip 1.6cm everywhere with equality somewhere at $t=0$;

\vskip 0.1cm\noindent \textbf{Type II}: The solution exists for
$t\in(-\infty,+\infty)$  and
$$|Rm|\leq 1 \ \ \ \ \ \ \ \ \
$$
\hskip 1.6cm everywhere with equality somewhere at $t=0$;

\vskip 0.1cm\noindent \textbf{Type III}: The solution exists for
$t\in(-A,+\infty)$ for some constant $A$ with $0<A<+\infty$ and
$$\ \ |Rm|\leq A/(A+t)
$$
\hskip 1.6cm everywhere with equality somewhere at $t=0$.
\end{definition}

\begin{definition} A solution of the Ricci flow is said to satisfy the
{\it injectivity radius condition} if for every sequence of {\it
(almost) maximum points} $\{(x_k,t_k)\}$, there exists a constant
$c_2>0$ independent of $k$ such that
$$inj(M,x_k,g_{ij}(t_k))\geq \frac{c_2}{\sqrt{K_{max}(t_k)}},
\quad \text{for all}\   k.$$ Here, by a sequence of  (almost)
maximum points, we mean $\{(x_k,t_k) \in M\times [0,T) \}$, $k=1,
2, \cdots$, has the following property: there exist positive
constants $c_1$ and $\alpha \in (0,1]$ such that
$$|Rm(x_k,t_k)|\geq c_1 K_{max}(t), \quad
t\in [t_k-\frac{\alpha}{K_{max}(t_k)},t_k]
$$ for all $k$.
\end{definition}

In \cite{P1}, Perelman proved an important no local collapsing
theorem, which yields the following result conjectured by Hamilton
in \cite{Ha95F}.

\begin{theorem} {\bf (Little Loop Lemma)}  Let $g_{ij}(t)$,
$0\le t<T<+\infty$, be a solution of the Ricci flow on a compact
manifold $M^n$. Then there exists a constant $\delta>0$ having the
following property: if at a point $x_0\in M$ and a time $t_0\in
[0,T)$, $$ |Rm|(\cdot,t_0)\leq r^{-2} \quad \text{on}\
B_{t_0}(x_0,r)
$$ for some $r\leq \sqrt{T}$, then the injectivity radius of $M$ with
respect to the metric $g_{ij}(t_0)$ at $x_0$ is bounded from below
by
$$inj(M,x_0,g_{ij}(t_0))\geq \delta r.$$
\end{theorem}

Clearly by the above Little Loop Lemma a maximal solution on a
compact manifold with the maximal time $T<+\infty$ always
satisfies the injectivity radius condition. Also, by the
Gromoll-Meyer injectivity radius estimate \cite{GM}, a solution on
a complete noncompact manifold with positive sectional curvature
also satisfies the injectivity radius condition. We refer the
reader to \cite{CZ05} (Chapter 4) for more detailed discussions.

\begin{theorem} {\bf (Hamilton \cite{Ha95F})} For any maximal solution
to the Ricci flow which satisfies
the injectivity radius condition and is of Type I, II, or III,
there exists a sequence of dilations of the solution which
converges in $C_{loc}^{\infty}$ topology to a singularity model of
the corresponding Type.
\end{theorem}

\vskip 0.2cm In the case of manifolds with nonnegative curvature
operator, or K\"ahler metrics with nonnegative holomorphic
bisectional curvature, we can bound the Riemannian curvature $Rm$
by the scalar curvature $R$ up to a constant factor depending only
on the dimension. Then we can slightly modify the statements in
the previous theorem as follows

\begin{theorem} {\bf (Hamilton \cite{Ha95F})} For any complete maximal
solution to the Ricci flow with bounded
and nonnegative curvature operator on a Riemannian manifold, or on
a K\"ahler manifold with bounded and nonnegative holomorphic
bisectional curvature, there exists a sequence of dilations which
converges to a singular model.

\vskip 0.1cm \noindent For Type I solutions:\  the limit model
exists for $t\in (-\infty,\Omega)$ with $0<\Omega<+\infty$ and has
$$R\leq \Omega/(\Omega-t)
$$ everywhere with equality somewhere at $t=0$;

\vskip 0.1cm \noindent For Type II solutions:\ the limit model
exists for $t\in (-\infty, +\infty)$  and has
$$R\leq 1 $$ everywhere with equality somewhere at $t=0$;

\vskip 0.1cm \noindent For Type III solutions:\ the limit model
exists for $t\in(-A, +\infty)$ with $0<A<+\infty$ and has
$$R\leq A/(A+t)$$ everywhere with equality somewhere at $t=0$.
\end{theorem}

For Type II or Type III singularity models with nonnegative
curvature we have the following results.

\begin{theorem} {\bf (Hamilton \cite{Ha93E})} Any Type II singularity model of
the Ricci flow with nonnegative curvature operator and positive
Ricci curvature must be a steady Ricci soliton.
\end{theorem}

\begin{theorem} {\bf (Cao \cite{Cao97})}

(i) Any Type II singularity model on a K\"ahler manifold with
nonnegative holomorphic bisectional curvature and positive Ricci
curvature must be a steady K\"ahler-Ricci soliton;

(ii) Any Type III singularity model on a K\"ahler manifold with
nonnegative holomorphic bisectional curvature and positive Ricci
curvature must be a shrinking K\"ahler-Ricci soliton.
\end{theorem}

\begin{theorem} {\bf (Chen-Zhu \cite{CZ00})} Any Type III singularity model
of the Ricci flow with nonnegative curvature operator and positive
Ricci curvature must be a homothetically expanding Ricci soliton.
\end{theorem}

We remark that the basic idea in proving the above theorems is to
apply the Li-Yau-Hamilton estimates for the Ricci flow \cite{Ha93}
or the K\"ahler-Ricci flow \cite{Cao92}, and the strong maximum
principle type arguments.

On the other hand, by exploring Perelman's $\mu$-entropy (defined
by Eq. (2.3)), Sesum \cite{Se} studied compact Type I singularity
model and obtained the following

\begin{theorem} {\bf (Sesum \cite{Se})} Let $(M,g_{ij}(t))$ be a
compact Type I singularity model obtained as a rescaling limit of
a Type I maximal solution. Then $(M,g_{ij}(t))$ must be a gradient
shrinking Ricci soliton.
\end{theorem}

\begin{remark}
Recently Naber \cite{Na} showed that a suitable rescaling
limit of any Type I maximal solution is a gradient shrinking
soliton. However, it is not guaranteed that such a limit soliton
is non-flat. It remains an interesting question whether a
noncompact Type I singularity model is a gradient shrinking
soliton.
\end{remark}

\section{Reduced Distance, Reduced Volume and Gradient Shrinkers}

In Section 2, we saw compact shrinkers are critical points of the
$\mu$-energy, which is monotone under the Ricci flow. In
\cite{P1}, Perelman also introduced a space-time distance function
$l(x, \tau)$, the reduced distance, which is analogous to the
distance function introduced by Li-Yau \cite{LY}. Perelman used
the reduced distance and its associated reduced volume to deal
with complete noncompact solutions with bounded curvature. He
established the comparison geometry to the reduced distance
function, and proved that the reduced volume is monotone under the
Ricci flow. This monotonicity formula is more useful for local
considerations, and is used to prove the noncollapsing theorem for
complete solutions to the Ricci flow with bounded curvature on
noncompact manifolds.

\subsection{The reduced distance}

We will write the Ricci flow in the backward version
$$
\frac{\partial}{\partial \tau}g_{ij}=2R_{ij}
$$
on a manifold $M$. In practice one often takes $\tau=t_0-t$ for
some fixed time $t_0$. Throughout the section we assume that $(M,
g_{ij}(\tau))$ is complete with bounded curvature\footnote{For
most parts of the discussion in this section, it suffices to
assume Ricci curvature bounded from below, see Ye \cite{Ye}.}.

To each space curve $\gamma(\tau)$ in $M$, $0\leq
\tau_1\le\tau\le\tau_2$, its {\bf $\mathcal{L}$-length}
 is defined as
$$
\mathcal{L}(\gamma) =\int_{\tau_1}^{\tau_2}\sqrt{\tau}\
[R(\gamma(\tau),\tau)+ |\dot{\gamma}(\tau)|^2_{g(\tau)}]\ d\tau.
$$
A curve $\gamma(\tau)$ is called an {\bf
$\mathcal{L}$-geodesic}\index{$\mathcal{L}$-geodesic} if the
tangent vector field $X=\dot{\gamma}(\tau)$ along $\gamma$
satisfies the {\bf $\mathcal{L}$-geodesic equation}
$$
\nabla_XX-\frac{1}{2}\nabla R+\frac{1}{2\tau}X+2 Rc(X,\cdot)=0.
$$

Given any space time point $(p, \tau_1)$ ($\tau_1\geq 0$) and a
tangent vector $v \in T_{p}M$, there exists a unique
$\mathcal{L}$-geodesic $\gamma(\tau)$ starting at $p$ with
$\lim\limits_{\tau\rightarrow
\tau_1}\sqrt{\tau}\dot{\gamma}(\tau)=v$. Also, for any $q\in M$
and $\tau_2>\tau_1$, there always exists an shortest
$\mathcal{L}$-geodesic $\gamma(\tau)$: $[\tau_1,\tau_2]\rightarrow
M$ connecting $p$ and $q$.

Now we fix a point $p\in M$ and set $\tau_1=0$. The {\bf
$\mathcal{L}$-distance function} on the space-time $M\times
\mathbb{R}^+$ is denoted by $L(q,\bar{\tau})$ and defined to be
the $\mathcal{L}$-length of the $\mathcal{L}$-shortest curve
$\gamma(\tau)$ connecting $p=\gamma (0)$ and
$q=\gamma(\bar{\tau})$.

The {\bf reduced distance} $l(q, \bar{\tau})$, from the space-time
origin $(p, 0)$ to $(q, \bar{\tau})$, is defined as

$$l(q, \bar{\tau})\triangleq \frac{1}{2 \sqrt{\bar{\tau}}}\; L(q,\bar{\tau}).$$

\begin{remark} In the case that $(M^n, g)$ is a static solution to the Ricci flow (i.e., $g$
is Ricci flat), it is easy to see that
$$L(q,\bar{\tau})=\frac{d^2(p,q)} {2\sqrt{\bar \tau}},$$
so $$ l(q,\bar{\tau})=\frac{d^2(p,q)} {4\bar \tau}.$$ In
particular, the reduced distance function $l(x, \tau)$ on $\Bbb
R^n$ with respect to the origin is given by the Gaussian shrinker
potential function (see Example 1.11(a))
$$ l(x,{\tau})=\frac{|x|^2} {4\tau}=\frac{|x|^2} {4(1-t)}.$$
\end{remark}

\begin{remark} One has $l(x, \tau)=n/2$ for positive Einstein manifold $(M^n, g(\tau))$,
normalized with $R=n/2\tau$, $\tau=1-t$. Furthermore, if $(M^n,
g_{ij}(x, t), f(x,t))$, $-\infty <t < 1$, is a complete gradient
shrinker with bounded curvature, then the reduced distance
function only differs from the potential function $f$ by a
constant: $l(x,\tau)=f(x, 1-\tau)+C$ with $\tau=1-t>0$, see, e.g.,
Lemma 7.77 in \cite{Cetc2}.
\end{remark}

In \cite{P1}, Perelman computed the first and second variations of
the $\mathcal{L}$-length and obtained

\begin{lemma}  \label{el}For the reduced distance $l(q,\bar{\tau})$
defined above, there hold

$$ \label{time} \frac{\partial l}{\partial
  \bar{\tau}}=-\frac{l}{\bar{\tau}}+R+\frac{1}{2\bar{\tau}^{3/2}}K
\eqno(4.1)$$
$$ \label{1space} |\nabla
l|^2=-R+\frac{l}{\bar{\tau}}-\frac{1}{\bar{\tau}^{3/2}}K
\eqno(4.2)$$
$$\label{2space}\Delta
  l\leq-R+\frac{n}{2\bar{\tau}}-\frac{1}{2\bar{\tau}^{3/2}}K. \eqno(4.3)$$
  where $$K=\int_0^{\bar{\tau}}\tau^\frac{3}{2}Q(X)d\tau,$$ and
   $$Q(X)=-R_\tau-\frac{R}{\tau}-2<\nabla R,X>+2{Rc}(X,X)$$
is the trace Li-Yau-Hamilton quadratic.  Moreover, the equality
holds in Eq.(4.3) if and only if the solution along the
$\mathcal{L}$ minimal geodesic $\gamma$ satisfies the gradient
soliton equation $$ R_{ij}+\nabla_i\nabla_j\ l-
\frac{1}{2\bar{\tau}}g_{ij}=0.
$$
\end{lemma}

\subsection {The reduced volume}

The reduced volume of $(M, g(\tau))$ is defined as

$$\tilde{V}(\tau)=\int_M (4\pi\tau)^{-\frac{n}{2}}\exp(-l(q,\tau))dvol_{\tau}(q).$$

\begin{remark} For the Gaussian soliton on $\Bbb R^n$, one finds

$$\tilde{V}(\tau)=1$$
for all $\tau>0$. Also, for a Ricci flat manifold $(M^n, g)$ one
has

$$\tilde{V}(\tau)=\int_M (4\pi\tau)^{-\frac{n}{2}}\exp(-d^2(p,q)/4\tau). \eqno(4.4)$$
\end{remark}

Note that it follows from Eqs. (4.1)-(4.3) in Lemma 4.1 that
$$
\label{eel} (\frac{\partial}{\partial
\tau}-\triangle+R)((4\pi\bar{\tau})^{-\frac{n}{2}}e^{-l})\leq 0.$$
Therefore $ \label{mone oe} \frac{d}{d
\bar{\tau}}\tilde{V}(\bar{\tau})\leq 0 $, provided $M$ is compact,
and the equality holds if and only if we are on a gradient
shrinking soliton.

More generally, Perelman \cite{P1} showed that the monotonicity of
$\tilde{V}(\bar{\tau})$ also holds on noncompact manifolds by
establishing a Jacobian comparison for the
$\mathcal{L}$-exponential map
$\mathcal{L}exp(\bar{\tau}):T_pM\rightarrow M$ associated to the
$\mathcal{L}$-length.

\begin{theorem} {\bf (Monotonicity of the reduced volume)}
Let $(M^n, g_{ij}(\tau))$ be a complete solution to the backward
Ricci flow $\frac{\partial}{\partial \tau}g_{ij}=2R_{ij}$ with
bounded curvature. Fix a point $p\in M$ and let $l(q,\tau)$ be the
reduced distance from $(p,0)$. Then

(i) the reduced volume $\tilde{V}(\tau)$
is nonincreasing in $\tau$, and $\tilde{V}(\tau)\leq 1$ for all
$\tau$;

(ii) the monotonicity is strict unless we are on a gradient
shrinking soliton.
\end{theorem}

\begin{remark}
As pointed out by Yokota \cite{Yo}, $\tilde{V}(\bar \tau)=1$ for
some $\bar \tau >0$ if and only if $(M^n, g_{ij}(\tau))$ is the
Gaussian shrinker on $\Bbb R^n$.
\end{remark}

Using the monotonicity of the reduced volume, one can derive a
version of no local collapsing theorem for noncompact solutions
with bounded curvature.

\begin{definition} Given any positive constants $\kappa>0$
and $r>0$, we say a solution to the Ricci flow is
$\kappa$-noncollapsed at $(x_0,t_0)$ on the scale $r$ if it
satisfies the following property: if $|Rm|(x,t)\leq r^{-2}$ for
all $(x, t)\in B_{t_0}(x_0,r)\times [t_0-r^{2},t_0],$ then
$$
vol_{t_0}(B_{t_0}(x_0,r))\geq \kappa r^{n}.
$$
\end{definition}

\begin{theorem} (Perelman \cite{P1}) \label{noncoll}
Let $(M^{n},\hat g_{ij})$ be a complete Riemannian manifold with
bounded curvature $|Rm|\leq k_0$ and with injectivity radius
bounded from below by $inj(M, \hat g_{ij})\geq i_0$, for some
positive constants $k_0$ and $i_0$. Let $g_{ij}(t),$ $t\in [0,T)$
be a smooth solution to the Ricci flow with bounded curvature for
each $t\in [0,T)$ and $g_{ij}(0)=\hat g_{ij}.$ Then there is a
constant $\kappa=\kappa (k_0, i_0, T)>0$ such that the solution is
$\kappa$-noncollapsed on scales $\leq \sqrt{T}.$
\end{theorem}

In the special case of gradient shrinking solitons, the above
result has been improved recently by Naber \cite{Na} and further
by Carrillo-Ni \cite{CN}.

\begin{proposition} {\bf (Carrillo-Ni \cite{CN})} Let $(M^n, g, f)$ be
a gradient shrinking soliton satisfying Eq. (1.2) for $\rho=1/2$
and with either bounded Ricci curvature $|Rc|\leq C$ or
nonnegative Ricci curvature $Rc\geq 0$. Then there exists a
positive constant $\kappa>0$ such that $vol(B(x_0,1))\geq \kappa$
whenever $|Rc|\leq 1$ on $B(x_0,1)$.
\end{proposition}

\begin{remark} In \cite{Na}, Naber showed a similar result but
requires the bound on the curvature tensor $Rm$.
\end{remark}

\section{Geometry of Gradient Ricci Solitons}

In this section, we examine the geometric structures of gradient
steady and shrinking Ricci solitons, in particular the sort we get
as Type I or Type II limits.

We start by examining ancient solutions with nonnegative
curvature, and two geometric quantities, the asymptotic scalar
curvature ratio and the asymptotic volume ratio.

\subsection{Ancient solutions of the Ricci Flow}

A complete solution $g_{ij}(t)$ to the Ricci flow is called {\it
ancient} if it is defined for $-\infty<t<T$. By definition, Type I
and Type II singularity models are ancient, and so are steady and
shrinking Ricci solitons.

Now let us recall the definitions of the asymptotic scalar
curvature ratio and the asymptotic volume ratio (cf.
\cite{Ha95F}).

Suppose $(M^n, g)$ is an $n$-dimensional complete noncompact
Riemannian manifold with nonnegative Ricci curvature $Rc\ge 0$.
Let $O$ be a fixed point on a Riemannian manifold $M^n$. We
consider the geodesic ball $B(O,r)$ centered at $O$ of radius $r$.
The well-known Bishop volume comparison theorem tells us that the
ratio $Vol(B(O,r))/r^n$ is nonincreasing in $r\in[0,+\infty)$.
Thus there exists a limit
$$\nu_M=\lim_{r\rightarrow+\infty}\frac{Vol(B(O,r))}{r^n},$$
which is called  {\it the asymptotic volume ratio} of $(M^n, g)$
which is invariant under dilation and is independent of the choice
of the origin.

\begin{remark} By Lemma 8.10 in \cite{Cetc2}, if $(M^n, g)$ is a
complete Ricci flat manifold, then the asymptotic volume ratio
$\nu_M$ is equal to the asymptotic reduced volume:
$$\nu_M=\lim_{\tau\to \infty}\tilde{V}(\tau).$$
\end{remark}

Next, let $s(x)$ denote the distance from $x\in M$ to the fixed
point $O$, and $R$ the scalar curvature. The {\it asymptotic
scalar curvature ratio} of $(M^n, g)$ is defined by
$$A=\limsup_{s\rightarrow +\infty}Rs^2,$$
which is also independent of the choice of the fixed point $O$ and
invariant under dilation.

\begin{remark}
The concept of asymptotic scalar curvature ratio is particular
useful on manifolds with positive sectional curvature. The first
gap type theorem was obtained by Mok-Siu-Yau \cite{MSY}.
 Yau (see \cite{GW}) suggested that this should be a general
 phenomenon, which confirmed by Greene-Wu \cite{GW}, Eschenberg-Shrader-Strake \cite{ESS} and
Drees \cite{Dr}. They showed that any complete noncompact
$n$-dimensional (except $n=4$ or 8) Riemannian manifold of
positive sectional curvature must have positive asymptotic scalar
curvature ratio $A>0$. Similar results on complete noncompact
K$\rm\ddot{a}$hler manifolds of positive holomorphic bisectional
curvature were obtained by Chen-Zhu \cite{CZ03} and Ni-Tam
\cite{NT}.
\end{remark}

\begin{proposition}{\bf (Hamilton \cite{Ha95F})} For a complete
ancient solution to the Ricci flow with bounded and nonnegative
curvature operator (or K\"ahler with bounded and nonnegative
bisectional curvature), the asymptotic scalar curvature ratio $A$
is constant (i.e., independent of time).
\end{proposition}

It is then natural to ask when is $A=\infty$, or $A<\infty$ for an
ancient solution. Intuitively, we see a paraboloid has $A=\infty$,
while a cone has $0<A<\infty$.

\begin{proposition} {\bf (Hamilton \cite{Ha95F})} Suppose $g_{ij}(t)$
is a complete ancient solution to the Ricci flow with bounded and
positive curvature operator. Assume $g_{ij}(t)$ is also Type-I
like so that
$$\sup_{t\in (-\infty, T)}\; R (T-t)<+\infty$$ and assume the
asymptotic scalar curvature ratio $A<\infty$. Then

(a) The asymptotic volume ratio $\nu_{M}$is positive; and

(b) for any origin $O$ and any time $t$, there exists $\phi (O,
t)>0$ such that $R(x, t)s(x)^2\ge \phi (O,t)$ at every point $x$,
where s(x) is the distance from $x$ to any origin $O$.
\end{proposition}

On the other hand, we have

\begin{theorem} {\bf (Perelman \cite{P1})}
Let $g_{ij}(t)$ be a complete non-flat ancient solution to the
Ricci flow with bounded and nonnegative curvature operator on a
noncompact manifold $M^n$. Then the asymptotic volume ratio
 $\nu_M(t)=0$ for all t.
\end{theorem}

In the K\"ahler case, the same result as Theorem 5.1 (assuming
nonnegative curvature operator) is obtained independently by
Chen-Zhu \cite{CZ04}. Moreover, we have the following stronger
result in which the assumption of nonnegativity of curvature
operator is replaced by the weaker assumption of nonnegative
holomorphic bisectional curvature.

\begin{theorem} {\bf (Chen-Tang-Zhu \cite{CTZ04} and Cao \cite{Cao04} for $n=2$,
Ni \cite{Ni} for $n\ge 3$)} Let $g_{\alpha\bar\beta}(t)$ be a
complete non-flat ancient solution to the K\"ahler-Ricci flow with
bounded and nonnegative holomorphic bisectional curvature on a
noncompact complex manifold $X^n$. Then the asymptotic volume
ratio $\nu_X(t)=0$ for all t.
\end{theorem}

Note that combining Theorem 5.1 with Proposition 5.2 we can deduce
the following

\begin{proposition}
Any complete noncompact ancient Type I-like solution to the Ricci
flow with bounded and positive curvature operator on an
$n$-dimensional manifold must have infinite asymptotic scalar
curvature ratio $A=\infty$.
\end{proposition}

\begin{remark}
Proposition 5.3 was first proved by Chow-Lu \cite{CL04} for $n=3$.
\end{remark}

Of course Theorem 5.1 implies that any complete non-flat shrinking
or steady soliton with bounded and nonnegative curvature operator
has zero asymptotic volume ratio $\nu_M=0$. On the other hand, for
gradient shrinking or expanding Ricci solitons, Carrillo-Ni
\cite{CN} recently showed the same result under weaker curvature
conditions.

\begin{proposition} {\bf (Carrillo-Ni \cite{CN})}
\medskip

(a) If $(M^n, g, f)$ is a gradient shrinking soliton with
nonnegative Ricci curvature $Rc\geq 0$, then the asymptotic volume
ratio $\nu_M=0$.

(b) If $(M^n, g, f)$ is a gradient expanding soliton with
nonnegative scalar curvature $R\geq 0$, then the asymptotic volume
ratio $\nu_M>0$.
\end{proposition}

\begin{remark}
It was known to Hamilton that if $(M^n, g, f)$ is a gradient
expanding soliton with bounded and nonnegative Ricci curvature
$0\leq Rc\leq C$, then the asymptotic volume ratio $\nu_M>0$.
\end{remark}

We like to point out that, due to the Hamilton-Ivey pinching
theorem\footnote{See, e.g, Theorem 2.4.1 in \cite{CZ05} or Theorem
6.44 in \cite{Cetc1})}, it turns out a $3$-dimensional complete
ancient solution $(M^3, g_{ij}(t))$ with bounded curvature
$|Rm|\le C$ has nonnegative sectional curvature. Also, every
complete ancient solution $(M^n, g_{ij}(t))$ with bounded
curvature $|Rm|\le C$ has nonnegative scalar curvature. This
follows from by applying the standard maximum principle to the
evolution equation of the scalar curvature $R$
$$\frac{\partial}{\partial t}R =\Delta R +2|Rc|^2.$$
Recently,  B.-L. Chen \cite{BChen} was able to remove the
curvature bound assumption in both results mentioned above.

\begin{theorem} {\bf (B.-L. Chen \cite{BChen})}
A $3$-dimensional complete ancient solution has nonnegative
sectional curvature.
\end{theorem}

\begin{proposition}\footnote{Though Proposition 5.5 was not stated
explicitly in \cite{BChen}, the essential arguments are there (see
Proposition 2.1(i) and Corollary 2.3(i) in \cite{BChen}).} {\bf
(B.-L. Chen \cite{BChen})} Let $g_{ij}(t)$ be a complete ancient
solution on a noncompact manifold $M^n$. Then the scalar curvature
$R$ of $g_{ij}(t)$ is nonnegative for all $t$.
\end{proposition}

To illustrate the local arguments used in \cite{BChen}, below we
sketch the proof of Proposition 5.5 which is relatively simpler.

\begin{proof} Suppose $g_{ij}$ is defined for $-\infty <t\leq T$
for some $T>0$. We can divide the arguments in \cite{BChen} into
two steps:

Step 1: Consider any complete solution $g_{ij}(t)$ defined on $[0,
T]$. For any fixed point $x_0\in M$, pick $r_0> 0$ sufficiently
small so that
$$|Rc|(\cdot, t)\leq (n-1) r_0^{-2} \qquad \mbox{on} \ B_t(x_0, r_0)$$ for
all $t\in [0, T]$. Then for any positive number $A>2$, pick
$K_A>0$ such that $R\geq -K_A$ on $B_0(x_0, Ar_0)$ at $t=0$. We
claim that there exists a universal constant $C>0$ (depending on
the dimension $n$) such that
$$R(\cdot, t)\geq \min\{-\frac{n} {t+\frac{1}{K_A}},
-\frac{C}{Ar_0^2}\} \qquad \mbox{on} \ B_t(x_0, \frac{3A}{4}r_0)
\eqno(5.1)$$ for each $t\in [0,T]$.

Indeed, take a smooth nonnegative decreasing function $\phi$ on
$\Bbb R$ such that $\phi =1$ on $(-\infty,7/8]$, and $\phi=0$ on
$[1, \infty)$. Consider the function $$u(x, t)=\phi(\frac{d_t(x_0,
x)} {Ar_0}) R(x, t).$$ Then we have
$$(\frac{\partial}{\partial t}-\Delta)u=\frac{\phi' R}{Ar_0}
(\frac{\partial}{\partial t}-\Delta)d_t(x_0,
x)-\frac{\phi''R}{(Ar_0)^2} +2\phi |Rc|^2-2\nabla \phi\cdot\nabla
R$$ at smooth points of the distance function $d_t(x_0, \cdot)$.

Let $u_{\min}(t)=\min_{M} u(\cdot, t)$. Whenever
$u_{\min}(t_0)\leq 0$, assume $u_{\min}(t_0)$ is achieved at some
point $\bar x\in B_{t_0}(x_0, r_0)$, then $\phi' R(\bar x,
t_0)\geq 0$. On the other hand, by Lemma 8.3(a) of Perelman
\cite{P1} (or Lemma 3.4.1(i), \cite{CZ05}), we know that
$$(\frac{\partial}{\partial t}-\Delta)d_t(x_0, x)\geq -\frac{5(n-1)}{3r_0}$$
outside $B_t(x_0, r_0)$. Following (Section 3, Hamilton
\cite{Ha86}), we define
$$\left.\frac{d}{dt}\right|_{t=t_0} u_{\min}=\liminf_{h\to 0^+}
\frac{u_{\min}(t_0+h)-u_{\min}(t_0)}{h},$$ the liminf of all
forward difference quotients, then

$$ \left.\frac{d}{dt}\right|_{t=t_0}
u_{\min}\geq -\frac{5(n-1)}{3Ar_0^2}\phi' R+ \frac{2}{n}\phi R^2 +
\frac{1}{(Ar_0)^2} (\frac{2{\phi'}^2}{\phi}-\phi'')R .$$ Hence,

$$\left.\frac{d}{dt}\right|_{t=t_0} u_{\min} \geq \frac{1}{n}u_{\min}^{2}(t_0)
-\frac{C^2}{(Ar_0^2)^2},$$ provided $u_{\min}(t_0)\leq 0$. Now
integrating the above inequality, we get
$$u_{\min}(t)\geq \min\{-\frac{n} {t+\frac{1}{K_A}},
-\frac{C}{Ar_0^2}\} \qquad \mbox{on} \ B_t(x_0, \frac{3A}{4}r_0),
$$ and the inequality (5.1) in our claim follows.

Step 2: Now if our solution $g_{ij}(t)$ is ancient, we can replace
$t$ by $t-\alpha$ in (5.1) and get
$$R(\cdot, t)\geq \min\{-\frac{n} {t-\alpha +\frac{1}{K_A}},
-\frac{C}{Ar_0^2}\} \qquad \mbox{on} \ B_t(x_0, \frac{3A}{4}r_0).
\eqno(5.2)$$ Letting $A\to \infty$ and then $\alpha \to -\infty$,
we completed the proof of Theorem 5.3.

\end{proof}

We also would like to mention Yokota \cite{Yo} recently proved an
interesting gap theorem for ancient solutions.

\begin{proposition}
There exists a (small) constant $\epsilon>0$ depending only on the
dimension $n$ which satisfies the following property: suppose
$(M^n, g(t))$, $-\infty<t\leq T_0$ is a complete ancient solution
to the Ricci flow with Ricci curvature bounded below such that
$$\lim_{\tau\to \infty}\tilde{V}(\tau)\geq 1-\epsilon.$$ Then
$(M^n, g(t))$ is a Gaussian shrinker on $\Bbb R^n$. Here
$\tilde{V}(\tau)$ is the reduced volume with respect to $(p, 0)$
for some base point $p\in M$, and $\tau=T_0-t$.
\end{proposition}

We end this subsection by noting the following classification
result of Hamilton \cite{Ha95F} on 2-dimensional ancient
$\kappa$-solutions.

\begin{theorem}

The only $2$-dimensional non-flat ancient solutions to the Ricci
flow with bounded and nonnegative curvature and
$\kappa$-noncollapsed on all scales are the round sphere and the
round real projective plane.
\end{theorem}

\subsection{Geometry of gradient steady and expanding Ricci solitons}

We now turn to gradient steady and expanding solitons.

\begin{proposition} {\bf (Hamilton \cite{Ha95F})} Suppose we have
a complete noncompact gradient steady  Ricci soliton $(M^n,
g_{ij})$ so that
$$R_{ij}=\nabla_i\nabla_jf$$
for some potential function $f$ on $M$. Assume the Ricci curvature
is positive $Rc > 0$, and the scalar curvature $R$ attains its
maximum $R_{max}$ at a point $x_0\in M^n$. Then
$$|\nabla f|^2+R=R_{max}$$
everywhere on $M^n$. Furthermore, the function $f$ is convex and
attains its minimum at $x_0$.
\end{proposition}

\begin{remark} It was observed by Hamilton and the author \cite{CH00}
that $f$ is also a exhaustion function of linear growth. As a
consequence, the underlying manifold of the gradient steady
soliton is diffeomorphic to the Euclidean space $\Bbb R^n$.
\end{remark}

\begin{remark}
In case of a complete gradient expanding soliton with nonnegative
Ricci curvature, the potential function $f$ is a convex exhaustion
function of quadratic growth. Hence we have

\begin{proposition}
Let $(M^n, g_{ij}, f)$ be a gradient expanding soliton with
$Rc\geq 0$. Then $M^n$ is diffeomorphic to $\Bbb R^n$.
\end{proposition}
\end{remark}
In the K\"ahler setting we have the following strong
uniformization type result.

\begin{proposition} {\bf (Bryant \cite{BR04} and Chau-Tam \cite{CT})}
Suppose we have a complete noncompact gradient steady
K\"ahler-Ricci soliton $(X^n, g_{i\bar j})$. Assume Ricci
curvature is positive $Rc>0$, and the scalar curvature $R$ attains
its maximum $R_{max}$ at a point $x_0\in X^n$. Then $X^n$ is
biholomorphic to the complex Euclidean space $\mathbb C^n$.
\end{proposition}

\begin{remark} Under the same assumption as in Proposition 5.9, it was
observed earlier by Hamilton and the author \cite{CH00} (see also
\cite{Cao04}) that $X^n$ is Stein and diffeomorphic to $\Bbb
R^{2n}$.
\end{remark}

\begin{remark}
Chau-Tam \cite{CT} also showed

\begin{proposition} Let $(X^n, g_{i\bar j})$ be a complete noncompact gradient
expanding K\"ahler-Ricci soliton with nonnegative Ricci curvature,
then $X^n$ is biholomorphic to $\mathbb C^n$.
\end{proposition}
\end{remark}

\begin{proposition} {\bf (Hamilton \cite{Ha95F})} For a complete
noncompact gradient steady Ricci soliton with bounded curvature
and positive sectional curvature of dimension $n\ge3$ where the
scalar curvature assume its maximum at a point $O\in M$, the
asymptotic scalar curvature ratio is infinite
$$A=\limsup_{s\rightarrow +\infty}Rs^2=+\infty.$$
\end{proposition}

\begin{remark} In fact, by examining Hamilton's proof in \cite{Ha95F}, one gets
the stronger conclusion that
$$A_{1+\epsilon}=\limsup_{s\rightarrow
+\infty}Rs^{1+\varepsilon}=+\infty$$ for arbitrarily small
$\varepsilon>0$
\end{remark}

\begin{remark} For $n=3$, also see the work of Chu \cite{Chu} for
more geometric information.
\end{remark}

One of the basic questions is to classify steady Ricci solitons
with positive curvature. When dimension $n=2$, we have the
following important uniqueness result of Hamilton.

\begin{theorem} {\bf (Hamilton \cite{Ha88})} The only
complete steady Ricci soliton on a two-dimensional manifold with
bounded (scalar) curvature $R$ which assumes its maximum
$R_{\max}=1$ at an origin is the cigar soliton on the plane
$\mathbb{R}^2$ with the metric
$$ds^2=\frac{dx^2+dy^2}{1+x^2+y^2}.$$
\end{theorem}

\begin{remark}
For $n=3$, Prelman (\cite{P1}, {\bf 11.9}) claimed that any
complete noncompact $\kappa$-noncollapsed gradient steady soliton
with bounded positive curvature must be the Bryant soliton. He
also conjectured that any complete noncompact three-dimensional
$\kappa$-noncollapsed ancient solution with bounded positive
curvature is necessarily a Bryant soliton.
\end{remark}

\subsection{Geometry of shrinking solitons}

In this subsection, we describe recent progress on complete
shrinking Ricic solitons.

First of all, by the splitting theorem of Hamilton \cite{Ha86}, a
complete shrinking Ricci soliton with bounded and nonnegative
curvature operator either has positive curvature operator
everywhere or its universal cover splits as a product $N\times
\mathbb R^{k}$, where $k\ge 1$ and $N$ is a shrinking soliton with
positive curvature operator. On the other hand, we know that
compact shrinking solitons with positive curvature operator are
isometric to finite quotients of round spheres, thanks to the
works of Hamilton \cite{Ha82, Ha86} (for $n=3, 4$) and
B\"ohm-Wilking \cite{BW} (for $n\ge 5$).

For dimension $n=3$, Perelman \cite{P1} proved the following

\begin{proposition} {\bf (Perelman \cite{P1})} There does not exist
a three-dimensional complete noncompact $\kappa$-noncollapsed
gradient shrinking soliton with bounded and positive sectional
curvature.
\end{proposition}

In other words, a three-dimensional complete $\kappa$-noncollapsed
gradient shrinking soliton with bounded and positive sectional
curvature must be compact.

Based on the above proposition, Perelman \cite{P2} obtained the
following important classification result (see also Lemma 6.4.1 in
\cite{CZ05}), which is an improvement of a result of Hamilton
(Theorem 26.5, \cite{Ha95F}).

\begin{theorem} {\bf (Perelman \cite{P1})} Let $g_{ij}(t)$ be a
nonflat gradient shrinking soliton to the Ricci flow on a
three-manifold $M^3$. Suppose $g_{ij}(t)$ has bounded and
nonnegative sectional curvature and is $\kappa$-noncollapsed on
all scales for some $\kappa>0$. Then $(M,g_{ij}(t))$ is one of the
following:

(i) the round three-sphere $\mathbb{S}^3$, or one of its metric
quotients;

(ii) the round infinite cylinder $\mathbb{S}^2\times \mathbb{R}$,
or its $\mathbb{Z}_2$ quotient.
\end{theorem}

Thus, the only three-dimensional complete noncompact
$\kappa$-noncollapsed gradient shrinking soliton with bounded and
nonnegative sectional curvature are either $\mathbb{R}^3$ or
quotients of $\mathbb{S}^2\times \mathbb{R}$.

In the past a few years, there has been a lot of attempts to
improve and generalize the above results of Perelman. Ni-Wallach
\cite{NW1} and Naber \cite{Na} dropped the assumption on
$\kappa$-noncollapsing condition and replaced nonnegative
sectional curvature by nonnegative Ricci curvature. In addition,
Ni-Wallach \cite{NW1} allows the curvature $|Rm|$ to grow as fast
as $e^{ar(x)}$, where $r(x)$ is the distance function to some
origin and $a>0$ is a suitable small positive constant. In
particular, they proved

\begin{proposition} {\bf (Ni-Wallach
\cite{NW1})} Any 3-dimensional complete noncompact non-flat
gradient shrinking soliton with nonnegative Ricci curvature
$Rc\geq 0$ and with $|Rm|(x)\leq Ce^{ar(x)}$ must be a quotient of
the round cylinder $\mathbb{S}^2\times \mathbb{R}$.
\end{proposition}

Recently, based on Theorem 5.3 and Proposition 5.13, B.-L. Chen,
X.-P. Zhu and the author \cite{CCZ08} observed that one can
actually remove all the curvature bound assumptions.

\begin{theorem} {\bf (Cao-Chen-Zhu \cite{CCZ08})}
Let $(M^3, g_{ij})$ be a $3$-dimensional complete noncompact
non-flat shrinking gradient soliton. Then $(M^3, g_{ij})$ is a
quotient of the round cylinder $\mathbb{S}^2\times \mathbb{R}.$
\end{theorem}

For $n=4$, Ni and Wallach \cite{NW2} showed that any
$4$-dimensional complete gradient shrinking soliton with
nonnegative curvature operator and positive isotropic curvature,
satisfying certain additional assumptions, is a quotient of either
$\mathbb{S}^4$ or $\mathbb{S}^3\times \mathbb{R}$. Partly based on
this result, Naber \cite {Na} proved

\begin{proposition} {\bf (Naber \cite{Na})}
Any 4-dimensional complete noncompact shrinking Ricci soliton with
bounded and nonnegative curvature operator is isometric to either
$\mathbb{R}^4$, or a finite quotient of $\mathbb{S}^3\times
\mathbb{R}$ or $\mathbb{S}^2\times \mathbb{R}^2$.
\end{proposition}

For higher dimensions,  H. Gu and X.-P. Zhu \cite {GZ} first
proved that any complete, \emph{rotationally symmetric}, non-flat,
$n$-dimensional ($n \geq 3$) shrinking Ricci soliton with
$\kappa$-noncollapsing on all scales and with bounded and
nonnegative sectional curvature must be the round sphere
$\mathbb{S}^n$ or the round cylinder $\mathbb{S}^{n-1} \times
\mathbb{R}$. Subsequently, Kotschwar \cite{KB} proved that the
only complete shrinking Ricci solitons (without the curvature sign
and bound assumptions) of rotationally symmetric metrics (on
$S^n$, $\mathbb{R}^{n}$ and $R\times S^{n-1}$) are, respectively,
the round, flat, and standard cylindrical metrics.

Recently, various authors proved classification results on
gradient shrinking solitons with vanishing Weyl curvature tensor
which include the rotationally symmetric ones as a special case.

\begin{proposition} {\bf (Ni-Wallach \cite{NW1})\footnote{See also \cite{CW08}
for an alternative proof.})} Let $(M^{n}, g)$ be a complete,
locally conformally flat gradient shrinking soliton with
nonnegative Ricci curvature. Assume that
$$|Rm|(x) \le e^{a(r(x) + 1)}$$ for some constant $a> 0$, where
$r(x)$ is the distance function to some origin. Then its universal
cover is  $\mathbb{R}^n$, $\mathbb{S}^n$ or
$\mathbb{S}^{n-1}\times \mathbb{R}$.
\end{proposition}

\begin{remark} In the compact case, Eminenti-La Nave-Mantegazza \cite{ELM}
first showed that every compact shrinking Ricci soliton with
vanishing Weyl tensor is a quotient of $\Bbb S^n$.
\end{remark}

In \cite{PW3}, Petersen-Wylie removed the assumption on
nonnegative Ricci curvature and replaced the pointwise curvature
growth assumption by certain integral bound on the norm of Ricci
tensor.

\begin{proposition} {\bf (Petersen-Wylie \cite{PW3})}
Let $(M^{n}, g)$ be a complete gradient shrinking Ricci soliton
with potential function $f$. Assume the Weyl tensor $W = 0$ and
$$\int_{M} |Rc|^{2} e^{-f} dvol<\infty,$$ then $(M^{n}, g)$ is a finite
quotient of $\mathbb{R}^n$, $\mathbb{S}^n$ or
$\mathbb{S}^{n-1}\times \mathbb{R}$.
\end{proposition}

On the other hand, Z.-H. Zhang \cite{Zh2} removed all curvature
bound assumptions in Proposition 5.15 and Proposition 5.16, thus
giving an extension of Theorem 5.7 to higher dimensions.

\begin{theorem} {\bf (Z.-H. Zhang \cite{Zh2})}
Any complete gradient shrinking soliton with vanishing Weyl tensor
must be a finite quotients of $\mathbb{R}^n$, $\mathbb{S}^n$ or
$\mathbb{S}^{n-1}\times \mathbb{R}$.
\end{theorem}

\begin{remark} In the K\"ahler case, we have a rather satisfying
result in all dimensions by Ni \cite{Ni}.

\begin{proposition} {\bf (Ni \cite{Ni})}
In any complex dimension, there is no complete noncompact gradient
shrinking K\"ahler-Ricci soliton with positive holomorphic
bisectional curvature.
\end{proposition}
\end{remark}

\begin{remark}
Wylie \cite{Wy} showed that a complete shrinking Ricci soliton has
finite fundamental group. In the compact case, the result was
first proved by Derdzinski \cite{De06} and Fern\'adez-L\'opez and
Garc\'ia-R\'io \cite{FG} (see also a different proof by
Eminenti-La Nave-Mantegazza \cite{ELM}). This can be viewed as an
extension of the well-known Myers' theorem. Moreover,
Fang-Man-Zhang \cite{FMZ} proved that a complete gradient
shrinking Ricci soliton with bounded scalar curvature has finite
topological type.
\end{remark}

\section{Stability of Ricci Solitons}

In this section we describe the second variation formulas for
Perelman's $\lambda$-energy and $\nu$-energy due to Hamilton,
Ilmanen and the author \cite{CHI04}.

\subsection{Second Variation of $\lambda$-energy} Recall that the
$\lambda$-energy is defined by
$$\lambda(g_{ij})=\inf\{{\mathcal{F}}(g_{ij},f): f\in
C^{\infty}(M), \int_Me^{-f}dV=1\}$$ and its first variation is
given by
\begin{align*}
\left.\frac{d}{ds}\right|_{s=0}\lambda(g(s))= \int_M
-h_{ij}(R_{ij}+\nabla_i\nabla_j f)e^{-f} dV,
\end{align*}
where $f$ is a minimizer.

For any symmetric 2-tensor $h=h_{ij}$ and 1-form
$\omega=\omega_i$, denote $Rm(h,h):=R_{ijkl}h_{ik}h_{jl}$,
$\text{div}\:\omega:=\nabla_i\omega_i$, $(\text{div}\:
h)_i:=\nabla_jh_{ji}$,
$(\text{div}^*\omega)_{ij}=-(\nabla_i\omega_j+\nabla_j\omega_i)/2
=-(1/2)L_{\omega^\#}g_{ij}$, where $\omega^\#$ is the vector field
dual to $\omega$.

\begin{theorem} {\bf (Cao-Hamilton-Ilmanen \cite{CHI04})} Let
$(M^n,g)$ be a compact Ricci flat manifold and consider variations
$g(s)=g+sh$. Then the second variation $\mathcal
D^2_g\lambda(h,h)$ of $\lambda$ at $g$ is given by
$$\left.\frac{d^2}{ds^2}\right|_{s=0}\lambda(g(s)) =\int_M <Lh,h> dV,$$ where
$$Lh:=\frac{1}{2}\Delta h+ \text{div}^*\: \text{div}\: h
+\frac{1}{2}{\nabla}^2v_{h}+Rm(h,\cdot),
$$
and $v_h$ satisfies
\begin{align*}
\Delta v_h=\text{div}\: \text{div}\: h.
\end{align*}
\end{theorem}

\noindent Note if we decompose $C^\infty(\text{Sym}^2(T^*M))$ as
$$
\ker\text{div}\oplus\text{im}\:\text{div}^*.
$$
One verifies that $L$ vanishes on $\text{im}\:\text{div}^*$, that
is, on Lie derivatives. On $\ker\text{div}$ one has
$$
L=\frac{1}{2}\Delta_L,
$$
where
$$
\Delta_Lh:=\Delta h+2Rm(h,\cdot)-Rc\cdot h-h\cdot Rc
$$
is the Lichnerowicz Laplacian on symmetric 2-tensors. We call a
critical point $g$ of $\lambda$ {\bf linearly stable} if $L\le0$.

\begin{example}
A Calabi-Yau K3 surface has $\Delta_L\le0$ according to
Guenther-Isenberg-Knopf \cite{GIK02}. More generally, Dai-Wang-Wei
\cite{DWW04} showed that any manifold with a parallel spinor has
$\Delta_L\le0$. So these manifolds are linearly stable in the
sense presented above.
\end{example}

\begin{example}
Let $g$ be compact and Ricci flat. Following \cite{B84, GIK02}, we
examine conformal variations. It is convenient to replace $ug$ by
$$
h=Su:=(\Delta u)g-D^2u
$$
which differs from the conformal direction only by a Lie
derivative and is divergence free. We have
$$
\Delta_LSu=(S\Delta u)g,
$$
so $\Delta_L$ has the same eigenvalues as $\Delta$. In particular,
$L\le 0$ in the conformal direction. This contrasts with the
Einstein functional.
\end{example}

\subsection{Second Variation of $\nu$-energy}

Recall that the $\nu$-energy is defined by $$
\nu(g_{ij})=\inf\{\mathcal W(g,f,\tau): f\in C^\infty(M), \tau>0,
\frac{1}{(4\pi\tau)^{n/2}}\int_M e^{-f}dV=1\}
$$ and its first variation is given by
\begin{align*}
\left.\frac{d}{ds}\right|_{s=0}\nu(g_{ij}(s))=\frac{1}{(4\pi\tau)^{n/2}}\int_M
-h_{ij}[\tau(R_{ij}+\nabla_i\nabla_j f)-g_{ij}/2]e^{-f} dV.
\end{align*}

We now describe the second variation of the $\nu$-energy in
\cite{CHI04} for positive Einstein metrics and compact shrinking
Ricci solitons.

First, for Einstein metrics normalized by $Rc=g/2\tau$, we have

\begin{theorem} {\bf (Cao-Hamlton-Ilmannen \cite{CHI04})} Let $(M,g)$ be a Einstein manifold
of positive scalar curvature and consider variations $g(s)=g+sh$.
Then the second variation $\mathcal D^2_g\nu(h,h)$ is given by
$$\left.\frac{d^2}{ds^2}\right|_{s=0}\nu(g(s))=\frac{\tau}{\text{vol}(g)}\int_M <Nh,h>,$$
where
$$
Nh:=\frac{1}{2}\Delta h+\mbox{div}^*\: \mbox{div}\: h+
\frac{1}{2}{\nabla}^2v_h+Rm(h,\cdot)-\frac{g}{2n\tau\mbox{vol}(g)}\int_M
\text{tr}_gh.
$$
and $v_h$ is the unique solution of
\begin{align*}
\Delta v_h+\frac{v_h}{2\tau}=\text{div}\: \text{div}\: h,
\qquad\int_M v_h=0.
\end{align*}
\end{theorem}

As in the previous case, $N$ is degenerate negative elliptic and
vanishes on $\im\div^*$.  Write
$$
\ker\div=(\ker\div)_0\oplus\RR g
$$
where $(\ker\div)_0$ is defined by $\int\tr_gh=0$. Then on
$(\ker\div)_0$ we have
$$
N=\frac{1}{2}\left(\Delta_L+\frac{1}{\tau}\right)
$$
where $\Delta_L$ is the Lichnerowicz Laplacian. So the linear
stability of a shrinker comes down to the (divergence free)
eigenvalues of the Lichnerowicz Laplacian. Let us write $\mu_L$
for the maximum eigenvalue of $\Delta_L$ on symmetric 2-tensors
and $\mu_N$ for the maximum eigenvalue of $N$ on $(\ker\div)_0$.
We now quote several examples from \cite{CHI04}:

\begin{example} The round sphere is linearly stable:
$\mu_N=-2/(n-1)\tau<0$. In fact, it is {\it geometrically stable}
(i.e. nearby metrics are attracted to it up to scale and gauge) by
the results of Hamilton \cite{Ha82,Ha86,Ha88} and Huisken
\cite{Hu}.
\end{example}

\begin{example}
\label{CPn-example} For complex projective space $\CP^m$, the
maximum eigenvalue of $\Delta_L$ on $(\ker\div)_0$ is
$\mu_L=-1/\tau$ by work of Goldschmidt \cite{G04}, so $\CP^m$ is
neutrally linearly stable, i.e. the maximum eigenvalue of $m$ on
$(\ker\div)_0$ is $\mu_N=0$.
\end{example}

\begin{example}
Any product of two Einstein manifolds $M=M_1^{n_1}\times
M_2^{n_2}$ is linearly unstable, with $\mu_N=1/2\tau$. The
destabilizing direction $h=g_1/n_1-g_2/n_2$ corresponds to a
growing discrepancy in the size of the factors.
\end{example}

\begin{example}
Any compact K\"ahler-Einstein manifold $X^n$ of positive scalar
curvature with $\dim H^{1,1}(X)\ge 2$ is linearly unstable.
Indeed, we can compute $\mu_N$ as follows. Let $\sigma$ be a
harmonic 2-form and $h$ be the corresponding metric perturbation;
then $\Delta_Lh=0$, and if $\sigma$ is chosen perpendicular to the
K\"ahler form, then as above we obtain $\mu_N=1/2\tau$.
\end{example}

\begin{example}
Let $Q^m$ denote the complex hyperquadric in $\CP^{m+1}$ defined
by
$$
\sum_{i=0}^{m+1}z_i^2=0,
$$
a Hermitian symmetric space of compact type, hence a
K\"ahler-Einstein manifold of positive scalar curvature.

(a) $Q^2$ is isometric to $\CP^1\times\CP^1$, the simplest example
of the above instability phenomenon.

(b) $Q^3$ has $\dim H^{1,1}(Q^3)=1$, so the above discussion does
not apply. But the maximum eigenvalue of $\Delta_L$ on
$(\ker\div)_0$ is $\mu_L=-2/3\tau$ by work of Gasqui and
Goldschmidt \cite{GG96} (or see \cite{GG04}). The proximate cause
is a representation that appears in the sections of the symmetric
tensors but not in scalars or vectors. Therefore, $Q^3$ is
linearly unstable with
$$
\mu_N=\frac{1}{6\tau}.
$$

(c) For $Q^4$, the maximum eigenvalue of $\Delta_L$ on symmetric
tensors is $\mu_L=-1/\tau$ by work of Gasqui and Goldschmidt
\cite{GG91} (or see \cite{GG04}). So $Q^4$ is neutrally linearly
stable: $\mu_N=0$.
\end{example}

For compact shrinkers, we have the following general second
variation formula.
\begin{theorem} {\bf (Cao-Hamilton-Ilmanen \cite{CHI04})} Let $(M,g)$ be a
shrinker with the potential function $F$. The second variation
$\DDD^2_g\nu(h,h)$ is given by
\begin{align*}
\left.\frac{d^2}{ds^2}\right|_{s=0}\nu(g(s))
&=\frac{\tau}{(4\pi\tau)^{n/2}}\int_M e^{-F} <\hat{N}h, h>,&\\
\end{align*}
where
$$
\hat Nh:=\frac{1}{2}\Delta_F h+\div_F^*\div_F h+
\frac{1}{2}D_F^2\hat
v_h+Rm(h,\cdot)-\frac{Rc}{n(4\pi\tau)^{n/2}}\int_M e^{-F} tr_gh.
$$
and $\hat v_h$ is the unique solution of
\begin{align*}
\Delta_F \hat v_h+\frac{\hat v_h}{2\tau}=\div_F\div_F
h,\qquad\int_M \hat v_h e^{-F}=0.
\end{align*}
\end{theorem}
Here we define
$$\div_F h=(e^{-F}\div e^F)h=(\div -D_jF)h=D_jh_{ij}-D_jF h_{ij},$$
$$\int_M e^{-F}(\phi, \div_F h)=\int_M e^{-F} (\div_F^*, h),$$
and
$$\Delta_F = \Delta -D_pF D_p.$$
It is easy to see that $\div_F^*=\div^*$.

Finally, my student Meng Zhu \cite{MengZhu} has computed the
corresponding second variation of the $\nu_{-}$-entropy for
compact negative Einstein manifolds.

\begin{theorem} Let $(M^n,g)$ be a compact Einstein manifold
of negative scalar curvature, normalized by $Rc=-g/2\sigma$. Then
the second variation $\mathcal D^2_g\nu_{-}(h,h)$ is given by
$$\left.\frac{d^2}{ds^2}\right|_{s=0}\nu_{-}(g(s))=\frac{\sigma}{\text{vol}(g)}\int_M <N_{-}h,h>,$$
where
$$
N_{-}h:=\frac{1}{2}\Delta h+\mbox{div}^*\: \mbox{div}\: h+
\frac{1}{2}{\nabla}^2v_h+Rm(h,\cdot)+\frac{g}{2n\sigma\mbox{vol}(g)}\int_M
\text{tr}_gh.
$$
and $v_h$ is the unique solution of
\begin{align*}
\Delta v_h-\frac{v_h}{2\sigma}=\text{div}\: (\text{div}\: h),
\qquad\int_M v_h=0.
\end{align*}
\end{theorem}

\begin{remark} There have been other works on related stability problems,
see, e.g., Chau \cite{Chau}, Chau-Schn\"urer \cite{CS},
Guenther-Isenberg-Knopf \cite{GIK06}, Sesum \cite{Se06},
Dai-Wang-Wei \cite{DWW07}, and Suneeta \cite{Su}.
\end{remark}

\section{The Gaussian Density of Shrinking Ricci Solitons}

In \cite{CHI04}, the notion of Gaussian density (or central
density) of a shrinking Ricci soliton is introduced. In case of a
compact shrinking soliton $(M^n, g)$, the Gaussian density is
simply given by
$$\Theta (M)=\Theta (M, g):=e^{\nu(M, g)}$$
where $\nu(M,g)=\nu(g_{ij})$ is the $\nu$-energy of the shrinker
$(M, g)$.

In the following discussion, we will normalize Einstein manifolds
of positive scalar curvature by $Rc=g/(n-1)$ (i.e., $\tau=\frac{1}
{2(n-1)}$), so that the round sphere $S^n$ has radius $1$. As
shown in \cite{CHI04}, we have the following facts.

(1) $\Theta(S^n)=\left(\dfrac{n-1}{2\pi e}\right)^{n/2}\vol(S^n)$.

(2) If $M$ is a Einstein manifold of positive scalar curvature,
then
$$
\Theta(M)=\left(\frac{1}{4\pi\tau e}\right)^{n/2}\vol_\tau(M)\le
\Theta(S^n),
$$
with equality if and only if $M=S^n$.

(3) $\Theta(\CP^m)=\left(\dfrac{m+1}{\pi
e}\right)^N\dfrac{\vol(S^{2m+1})}{2\pi}$.

(4) The K\"ahler-Einstein manifold $M=\CP^2\#k(-\CP^2)$,
$k=0,3,\ldots,8$,  has $\Theta(M)=(9-k)/2e^2$.

(5) $\Theta(M_1\times M_2)=\Theta(M_1)\Theta(M_2)$.

\vskip 0.1cm As in \cite{CHI04}, we say that one shrinking soliton
{\it decays} to another if there is a small perturbation of the
first whose Ricci flow develops a singularity modelled on the
second.  Because the $\nu$-invariant is monotone during the flow,
decay can only occur from a shrinking soliton of lower density to
one of higher density. This creates a ``decay lowerarchy''.

\section{$4$-D Einstein Manifolds and Shrinking Ricci Solitons} In this
section we collect information about stability and Gaussian
density values of all known (orientable) positive Einstein
$4$-manifolds and $4$-dimensional compact shrinking Ricci
solitons. Below is a list containing all the information we know
so far. Note that on $\CP^2\#(-\CP^2)$ we have the K\"ahler-Ricci
soliton metric of \cite{Ko, Cao94} and the Page non-K\"ahler
Einstein metric \cite{Pa78}, while on $\CP^2\#2(-\CP^2)$ we have
the K\"ahler-Ricci soliton metric of \cite{WZ} and the very recent
Non-K\"ahler Einstein metric of \cite{CLW}. Also, Headrick-Wiseman
\cite{HW} found out $\Theta=.455$ for the Wang-Zhu soliton through
numerical computation\footnote{Added in Proof: Very recently
Stuart Hall has found the value of $\Theta=.4552$ for the
Non-K\"ahler Einstein metric in \cite{CLW}. His calculation also
yielded the value of $\Theta=.4549$ for Wang-Zhu soliton, in line
with \cite{HW}.}.

\begin{remark} On any compact Einstein $4$-manifold $M^4$, the Hitchin-Thorpe
inequality (see e.g., \cite{Besse}) states that
$$2\chi(M)\geq 3|\tau(M)|,$$
where $\chi(M)$ is the Euler characteristic  and $\tau(M)$ is the
signature of $M^4$. An interesting question is whether the
Hitchin-Thorpe inequality holds for $4$-dimensional compact
shrinkers. Note that for the K\"ahler-Ricci shrinkers
$M=\CP^2\#k(-\CP^2)$ ($k=1, 2$), $\chi(M)=3+k$ and
$|\tau(M)|=k-1$, so $2\chi(M)\geq 3|\tau(M)|$ is valid.
\end{remark}

\bigskip
\begin{tabular}
[c]{|llllr|} \hline
&&&&\\
\textbf{Shrinking Solitons} & \textbf{Type} & $\Theta$ & $\Theta$ &\textbf{Stability}\\
\hline
&&&&\\
$S^4$ & Einstein & ${6/e^2}$ & .812 & Stable \\    
&&&&\\
$\CP^2$ & Einstein & ${9/2e^2}$ & .609 & Stable \\  
&&&&\\
$S^2\times S^2$ & Einstein product & ${4/ e^2}$ & .541 & Unstable \\ 
&&&&\\
$\CP^2\#(-\CP^2)$ & K\"ahler-Ricci soltion (\cite{Ko},\cite{Cao94}) & ${3.826/e^2}$ & .518 & Unknown \\
&&&&\\
$\CP^2\#(-\CP^2)$ & non-K\"ahler Einstein  (\cite{Pa78}) & ${3.821/e^2}$ & .517 & Unknown \\
&&&&\\
$\CP^2\#2{-\CP^2}$ & K\"ahler-Ricci soltion (\cite{WZ}) & ${3.361/e^2}$ & .4549 & Unknown \\
&&&\\
$\CP^2\#2{-\CP^2}$ & non-K\"ahler Einstein (\cite{CLW}) &
${<7/2e^2}$ & .4552 & Unknown \\
&&&&\\
$\CP^2\#3(-\CP^2)$ & K\"ahler-Einstein & ${3/e^2}$ & .406 & Unstable \\  
&&&&\\
$\CP^2\#4(-\CP^2)$ & K\"ahler-Einstein & ${5/2e^2}$ & .338 & Unstable \\  
&&&&\\
$\CP^2\#5(-\CP^2)$ & K\"ahler-Einstein & ${2/e^2}$ & .271 & Unstable \\  
&&&&\\
$\CP^2\#6(-\CP^2)$ & K\"ahler-Einstein & ${3/2e^2}$ & .203 & Unstable \\  
&&&&\\
$\CP^2\#7(-\CP^2)$ & K\"ahler-Einstein & ${1/e^2}$ & .135 & Unstable \\  
&&&&\\
$\CP^2\#8(-\CP^2)$ & K\"ahler-Einstein & ${1/2e^2}$ & .068 & Unstable \\  
&&&&\\
\hline
\end{tabular}
\bigskip

\section{Open Problems}

We conclude with several problems/questions:

\vskip 0.1cm (i) Find a non-product, non-Einstein, Riemannian
(non-K\"ahler) shrinking Ricci soliton.

\vskip 0.1cm (ii) Show that any noncompact Type I singularity
model obtained as a rescaling limit of a Type I maximal solution
is a gradient shrinking soliton. Naturally, one should try to
explore Perelman's reduced volume approach.

\vskip 0.1cm (iii) Show that any $3$-dimensional complete
noncompact gradient steady soliton with bounded and positive
curvature is the Bryant soliton. Are there non-rotationally
symmetric steady soliton with positive curvature in dimensions
$n\ge 4$?

\vskip 0.1cm (iv) For $n\ge 4$, are there any complete noncompact
($\kappa$-noncollapsed) gradient shrinking Ricci soliton with
positive curvature operator or positive Ricci curvature? So far
there seems no known example of noncompact non-product shrinker
with nonnegative Ricci curvature.

\vskip 0.1cm (v) Is it true that linearly stable compact
$4$-dimensional shrinking solitons are necessarily Einstein? Are
the only linearly stable Einstein $4$-manifolds either the round
sphere $S^4$ or the complex projective space $\CP^2$ with the
Fubini-Study metric?

\vskip 0.1cm (vi) Does the Hitchin-Thorpe inequality hold for
compact 4-dimensional shrinking solitons?

\end{document}